\documentclass[11pt,a4paper]{article}
\usepackage{bbm}
\usepackage{mathrsfs}
\usepackage{amsfonts}
\usepackage{}

\usepackage{graphicx}
\usepackage{amsfonts,amsmath, amssymb}
\usepackage{color}

\newtheorem{theo}{\textbf{\ \ \quad Theorem}}[section]

\newtheorem{lem}{\textbf{\ \ \quad Lemma}}[section]
\newtheorem{remark}{\textbf{\ \ \quad Remark}}[section]

\newcommand{\lbl}[1]{\label{#1}}

\newcommand{\be}{\begin{equation}}
\newcommand{\ee}{\end{equation}}
\newcommand\bes{\begin{eqnarray}}
\newcommand\ees{\end{eqnarray}}
\newcommand{\bess}{\begin{eqnarray*}}
\newcommand{\eess}{\end{eqnarray*}}

{\setlength\arraycolsep{2pt}}

\title{Schauder and Sobolev Estimates of Parabolic Equations}

\author{Guangying Lv$^a$,  Jinlong Wei$^b$\\
\\
{\small \it $^a$ Institute of Applied Mathematics, Henan University} \\ {\small \it Kaifeng, Henan 475001, China}\\ {\small \tt gylvmaths@henu.edu.cn}\\
   {\small \it $^b$ School of Statistics and Mathematics, Zhongnan University of}\\
   {\small \it Economics and Law, Wuhan 430073, China}\\
    {\small \tt  weijinlong.hust@gmail.com}}

\begin{document}
\maketitle

\medskip

\begin{abstract} In this note, we use the non-homogeneous Poisson stochastic
process to show how knowing Schauder and Sobolev  estimates for
the one-dimensional heat equation allows one to derive their
multidimensional analogs. The method is probability. We generalize the result of Krylov-Priola \cite{KP2018}.



\end{abstract}

\baselineskip=15pt

\section{Introduction}
\setcounter{equation}{0}

For the classical theory of partial differential equations, the
Schauder and Sobolev estimates are important issues, see the book
\cite{Cbook2003,Ebook1998}. In \cite{KP2018}, Krylov-Priola used the Poisson
stochastic process to obtain the Schauder and Sobolev estimates of multi-dimensional heat equation
from the one-dimensional case. More precisely,  they first know the
Schauder and Sobolev estimates for the following equation
  \bes \left\{\begin{array}{llll}
\partial_tu(t,x)=D^2_xu(t,x)+f(t,x),\ \ t\in(0,T), &x\in\mathbb{R},\\
u(0,x)=0, \ &x\in\mathbb{R},
  \end{array}\right.
   \lbl{1.1}\ees
then they derive the Schauder and Sobolev estimates for multi-dimensional equation.
Actually, they obtained more abundant results.

The regularity of partial differential equations has been studied by many authors,
steady transport equation \cite{GT2010}, stochastic evolution equations \cite{BH2016,ZLW2018}, stochastic 
partial differential equations \cite{DL2016,DL2019} and so on. There are a lot of work about regularity
focusing on stochastic process, for example \cite{L2012,LZ2016}.

In the present paper, we aim to use the non-homogeneous Poisson stochastic
process to find some new results. The main difference between this paper and \cite{KP2018} is that
we use the non-homogeneous Poisson stochastic
process but Krylov-Priola used the homogeneous Poisson stochastic
process. The method used in \cite{KP2018} is probability and the results are interesting.

Throughout this paper, $T$ is a fixed positive number, $\mathbb{R}^d$ denotes Euclidean space and
$C^\alpha(\mathbb{R}^d)$, $\alpha\in(0,1)$ is the space of all real-valued functions
$f$ on $\mathbb{R}^d$ with the norm
  \bess
\|f\|_{C^\alpha(\mathbb{R}^d)}=\sup_{x\in\mathbb{R}^d}|f(x)|+[f]_{C^\alpha(\mathbb{R}^d)}<+\infty,
   \eess
where
    \bess
[f]_{C^\alpha(\mathbb{R}^d)}=\sup_{x\neq y}\frac{|f(x)-f(y)|}{|x-y|^\alpha}.
   \eess
As usual, we denote $C^{2+\alpha}(\mathbb{R}^d)$ as the space of real-valued twice continuously differentiable
functions $f$ on $\mathbb{R}^d$ with the norm
  \bess
\|f\|_{C^{2+\alpha}(\mathbb{R}^d)}=\sup_{x\in\mathbb{R}^d}(|f(x)|+|Df(x)|+|D^2f(x)|)+[D^2f]_{C^\alpha(\mathbb{R}^d)},
   \eess
where $Df$ is the gradient of $f$ and $D^2f$ is its Hessian.

The rest of this paper is arranged as follows. In Sections 2, we present some preliminaries and main result. Section 3 is the proof of main result.

\section{Preliminaries and Main Results}\label{sec2}\setcounter{equation}{0}
Consider the following Cauchy problem
     \bes \left\{\begin{array}{llll}
\partial_tu(t,x)=a(t)D^2_xu(t,x)+f(t,x),\ \ t\in(0,T), &x\in\mathbb{R},\\
u(0,x)=0, \ &x\in\mathbb{R},
  \end{array}\right.
   \lbl{2.1}\ees
where $a(t)$ is a positive bounded function. Denote $B_c((0,T),C_0^\infty(\mathbb{R}^d))$ as the space of functions
$\varphi$ satisfying that $\varphi$ is Borel bounded function and $\varphi(t,\cdot)\in C_0^\infty(\mathbb{R}^d)$ for any
$t\in(0,T)$.

It follows from \cite{Ebook1998,Krylov2008,Lieberman1996} that if $f$ belongs to $B_c((0,T),C_0^\infty(\mathbb{R}))$,
then (\ref{2.1}) has a solution $u(t,x)$ satisfying

  (i) $u$ is a continuous function in $[0,T]\times\mathbb{R}$;

  (ii) for any fixed $t\in[0,T]$, $u$ belongs to $C^{2+\alpha}(\mathbb{R})$ and has the following estimate
    \bes
\sup_{t\in[0,T]}\|u(t,\cdot)\|_{C^{2+\alpha}(\mathbb{R})} \leq N(T,\alpha)\sup_{t\in[0,T]}\|f(t,\cdot)\|_{C^{\alpha}(\mathbb{R})}.
   \lbl{2.2}\ees
Moreover, there exists only one solution $u$ satisfying the following properties
   \bes
&&\sup_{(t,x)\in[0,T]\times\mathbb{R}}|u(t,x)|\leq T\sup_{(t,x)\in[0,T]\times\mathbb{R}}|f(t,x)|,\lbl{2.3}\\
&&\sup_{t\in[0,T]}[D^2_xu(t,\cdot)]_{C^{\alpha}(\mathbb{R})}\leq N(\alpha)\sup_{t\in[0,T]}[f(t,\cdot)]_{C^{\alpha}(\mathbb{R})},\lbl{2.4}\\
&&\|D^2_xu\|^p_{L^p((0,T)\times\mathbb{R})}\leq N_p\|f\|^p_{L^p((0,T)\times\mathbb{R})}.
  \lbl{2.5}\ees
Here $L^p$-space is defined as usual.

Now we recall some knowledge of Poisson  stochastic process.
A non-homogeneous Poisson process $\pi(t,\omega)$ ($\pi_t$ for short) is a Poisson process with rate parameter $\lambda(t)$
such that the rate parameter of the process is a function of time. The significant difference
between the homogeneous and non-homogeneous Poisson process is that the latter case is
not a stationary process. Thus we can not write the non-homogeneous Poisson process as
the sum of  a sequence which is an i.i.d
(independently identically distribution) random variables.

As usual, $\pi_t$ is a counting process with the following properties
   \bess
{\rm(i)}:\ \mathbb{P}(\pi_t-\pi_s=k)=\frac{[m(t)-m(s)]^k}{k!}e^{-[m(t)-m(s)]},\ m(t)=\int_0^t\lambda(s)ds; \qquad\qquad\qquad\quad
  \eess

(ii) $\pi_t-\pi_s$ is independent of the trajectory $\{\pi_r,r\in[0,s]\}$.

For simplicity, in this paper, we only consider the 2-dimensional heat equation. For $x,y\in \mathbb{R}$, we set $z=(x,y)\in\mathbb{R}^2$.
For $l\in\mathbb{R}^2$, denote $D^2_l=l^il^jD_{ij}$, $D_i=D_{x_i}=\partial/\partial{x_i}$ and
$D_{ij}=D_iD_j$, where $i,j=1,2$ and $x_1=x,x_2=y$.
We obtain the following result.
  \begin{theo}\lbl{t2.1}
Let $a(t)>0$ be a bounded Borel measurable function. Then for
any $f\in B_c((0,T),C_0^\infty(\mathbb{R}^2))$, there exists a unique continuous in
$[0,T]\times\mathbb{R}^2$ solution $v(t,z)$of the equation
  \bes\left\{\begin{array}{llll}
\partial_tv(t,z)=a(t)\Delta v(t,z)+f(t,z),\ t>0,\ &z\in\mathbb{R}^2,\\
v(0,z)=0,\ \ \  &z\in\mathbb{R}^2.
   \end{array}\right.
   \lbl{2.6} \ees
Moreover, $v(t,\cdot)\in C^{2+\alpha}(\mathbb{R}^2)$ satisfies
   \bess
&&\sup_{(t,z)\in[0,T]\times\mathbb{R}^2}|v(t,z)|\leq T\sup_{(t,z)\in[0,T]\times\mathbb{R}^2}|f(t,z)|,\\
&&\sup_{t\in[0,T]}[D_{ij}v(t,\cdot)]_{C^{\alpha}(\mathbb{R}^2)}\leq N_0(\alpha)\sup_{t\in[0,T]}[f(t,\cdot)]_{C^{\alpha}(\mathbb{R}^2)},\\
&&\sup_{(t,z)\in[0,T]\times\mathbb{R}^2}[D^2_{l}v(t,z+l\cdot)]_{C^{\alpha}(\mathbb{R}^2)}\leq
N_0(\alpha)\sup_{(t,z)\in(0,T)\times\mathbb{R}^2}[D^2_{l}f(t,z+l\cdot)]_{C^{\alpha}(\mathbb{R}^2)},\\
&&\|D^2_{l}u\|^p_{L^p((0,T)\times\mathbb{R}^2)}\leq N_p\|f\|^p_{L^p((0,T)\times\mathbb{R}^2)},
  \eess
where $N_0(\alpha)$ and $N_p$ are positive constants.
  \end{theo}

\begin{remark}\lbl{r2.1}
The result of this paper has a little difference from \cite{KP2018} in the following part.
If $a(t)=1$, that is, $\lambda(t)\equiv\lambda $, then Theorem \ref{t2.1} is exactly the second
part of \cite{KP2018}. The big difference is that we can assume $\lambda(t)=h^2a(t)$ and then the
equation will keep the same form as the dimensional case. Of course, in \cite[Section 3]{KP2018},
Krylov-Priola used a suitable transform to consider the problem (\ref{2.1}). Here we
emphasize that we can use another stochastic process to deal with the problem (\ref{2.1}).

One can use renew process to study the regularity of parabolic equations. The difference
is that in the following Lemma \ref{l3.1}, $\mathbb{E}[\pi_{(k+1)2^{-n}}-\pi_{k2^{-n}}]$ will be
different. But for parabolic equation, the Poisson process is the best choice.
\end{remark}

\section{The Proof of Theorem \ref{t2.1}}\label{sec3}\setcounter{equation}{0}

In this section, we prove the main result. Similar to \cite{KP2018}, we consider the following equations
  \bes\left\{\begin{array}{llll}
\partial_tu(t,x,y,\omega)=a(t)D_x^2u(t,x,y,\omega)+f(t,x,y-h\pi_t(\omega)),\ t>0,\ &x\in\mathbb{R},y\in\mathbb{R},\\
u(0,x,y)=0,\ \ \  &x\in\mathbb{R},y\in\mathbb{R},
   \end{array}\right.
   \lbl{3.1} \ees
where $a(t)>0$ is a bounded Borel measurable function and  $h\in\mathbb{R}$ is a parameter. As usual in probability theory,
we do not indicate the dependence on $\omega$ in the sequence.  From the result of one-dimensional case, we get that
there exists a unique solution $u(t,x,y)$, depending on $y$ and $\omega$ as parameters.
And thus estimates (\ref{2.2})-(\ref{2.5}) hold for each $\omega\in\Omega$ and $y\in\mathbb{R}$ if
we replace $u(t,x)$ and $f(t,x)$ with $u(t,x,y)$ and $f(t,x,y-h\pi_t)$, respectively.

The solution of (\ref{3.1}) can be written as
   \bes
u(t,x,y+h\pi_t)=\int_0^t[a(s)D_x^2u(s,x,y+h\pi_s)+f(s,x,y)]ds+\int_{(0,t]}g(s,x,y)d\pi_s,
   \lbl{3.2}\ees
where
   \bes
g(s,x,y)=u(s,x,y+h+h\pi_{s-})-u(s,x,y+h\pi_{s-})
    \lbl{3.3} \ees
is the jump of the process $u(t,x,y+h\pi_t)$ as a function of $t$ at moment $s$ if $\pi_t$ has
a jump at $s$. Here $\pi_{s-}$ is the left-continuous w.r.t. $s$.

In order to prove the main result, we need to study the function $g$.
  \begin{lem}\lbl{l3.1} For $g$ defined as (\ref{3.3}) and $t\leq T$ we have
     \bess
\mathbb{E}\int_{(0,t]}g(s,x,y)d\pi_s=\int_0^t\lambda(s)[v(s,x,y+h)-v(s,x,y)]ds,
   \eess
where
   \bess
v(t,x,y):=\mathbb{E}u(t,x,y+h\pi_t).
   \eess
 \end{lem}

{\bf Proof.} Assume that $t=1$ for simplicity. Fix $x$ and $y$, and denote $g(s)=g(s,x,y)$.
Note that $g$ is bounded on $\Omega\times(0,T)$, and thus if we  define
   \bess
g_n(s)=g(k2^{-n})=u(k2^{-n},x,y+h+h\pi_{k2^{-n}-})-u(k2^{-n},x,y+h\pi_{k2^{-n}-})
  \eess
for $s\in(k2^{-n},(k+1)2^{-n}]$, $k=0,1,\dots$, then $g_n(s)\to g(s)$ as $n\to\infty$
for any $s\in(0,t]$ and $\omega\in\Omega$, and
   \bess
\xi_n:=\int_{(0,t]}g_n(s)d\pi_s\to\int_{(0,t]}g(s)d\pi_s=:\xi
   \eess
for any $\omega\in\Omega$. Dominated convergence theorem implies that $\mathbb{E}\xi_n\to\mathbb{E}\xi$.

Notice that
  \bes
\mathbb{E}\xi_n=\sum_{k=0}^{2^n-1}\mathbb{E}g(k2^{-n})(\pi_{(k+1)2^{-n}}-\pi_{k2^{-n}}).
   \lbl{3.4}\ees
Since the non-homogeneous Poisson process is an independent increment process, the
expectations of he products on the right in (\ref{3.4}) are equal to the products of expectations,
and since $\mathbb{E}\pi_t=m(t)$, we arrive at
  \bess
\mathbb{E}\xi_n&=&\mathbb{E}\sum_{k=0}^{2^n-1}g(k2^{-n})[m(k+1)2^{-n}-m(k2^{-n})]
=\mathbb{E}\int_0^tg_n(s)\lambda(s)ds\\
&&\to \mathbb{E}\int_0^tg(s)\lambda(s)ds
=\int_0^t\lambda(s)\mathbb{E}g(s)ds.
  \eess
Noting that for any $s>0$, we have $\pi_s=\pi_{s-}$ almost surely, and thus
   \bess
\mathbb{E}g(s)=v(s,x,y+h)-v(s,x,y).
  \eess
The proof is complete. $\Box$

Taking expectations on both sides of (\ref{3.2}), we obtain the following result.

\begin{lem}\lbl{l3.2} Let $f\in B_c(0,T),C_0^\infty(\mathbb{R}^2)$, $h\in\mathbb{R}$
and $\lambda(t)>0$ for all $t\in[0,T]$. Then there exists a unique continuous function
$v(t,x,y)$, $t\in[0,T]$, $x,y\in\mathbb{R}$, satisfying the equation
\bes
\partial_tv(t,x,y)=a(t)D_x^2v(t,x,y)+\lambda(t)[v(t,x,y+h)-v(t,x,y)]+f(t,x,y)
   \lbl{3.5}\ees
for $t\in(0,T)$, $x,y\in\mathbb{R}$, with zero initial condition and such that
$v(t,\cdot,y)\in C^{2+\alpha}(\mathbb{R})$ for any $t\in(0,T)$, $y\in\mathbb{R}$ and
   \bess
\sup_{(t,y)\in[0,T]\times\mathbb{R}}\|v(t,\cdot,y)\|_{C^{2+\alpha}(\mathbb{R})}\leq
N(T,\alpha)\sup_{(t,y)\in[0,T]\times\mathbb{R}}\|f(t,\cdot,y)\|_{C^{\alpha}(\mathbb{R})}.
   \eess
Furthermore,
   \bess
&&\sup_{(t,z)\in[0,T]\times\mathbb{R}^2}|v(t,z)|\leq T\sup_{(t,z)\in[0,T]\times\mathbb{R}^2}|f(t,z)|,\lbl{2.3}\\
&&\sup_{(t,y)\in[0,T]\times\mathbb{R}}[D_x^2v(t,\cdot,y)]_{C^{\alpha}(\mathbb{R})}\leq N(\alpha)\sup_{(t,y)\in(0,T)\times\mathbb{R}}[f(t,\cdot,y)]_{C^{\alpha}(\mathbb{R})},\lbl{2.4}\\
&&\|D_x^2v\|^p_{L^p((0,T)\times\mathbb{R}^2)}\leq N_p\|f\|^p_{L^p((0,T)\times\mathbb{R}^2)}.
  \lbl{3.5}\eess
\end{lem}

The proof of this lemma is similar to \cite[Lemma 2.2]{KP2018} and we omit it here.

Next, we will do with (\ref{3.5}) almost the same thing as with (\ref{2.1}). More precisely, we
consider $v(t,x,y)$ depending on $\omega$ as a unique solution of
   \bess
\partial_tv(t,x,y)=a(t)D_x^2v(t,x,y)+\lambda(t)[v(t,x,y+h)-v(t,x,y)]+f(t,x,y+h\pi_t)
    \eess
with zero initial condition. Then it follows from the above computations, we have the
function $w(t,x,y)=\mathbb{E}v(t,x,y-h\pi_t)$ satisfies
  \bes
\partial_tw(t,x,y)=a(t)D_x^2w(t,x,y)+\lambda(t)[w(t,x,y+h)-2w(t,x,y)+w(t,x,y-h)]+f(t,x,y).
    \lbl{3.6}\ees
Furthermore, $w(t,x,y)$ has the same estimates as in Lemma \ref{l3.2}.

{\bf Proof of Theorem \ref{t2.1}} Taking $\lambda(t)=h^2a(t)$ in (\ref{3.6}) and
letting $h\to0$, we have the solution $w=w_h$ of (\ref{3.6}) will converge to a function
$v(t,x,y)$, which satisfies the equation (\ref{2.6}). Furthermore, $v$ is continuous in $[0,T]\times\mathbb{R}^2$, and
is infinitely differentiable w.r.t. $(x,y)$ for any $t\in(0,T)$ and all the estimates in Lemma
\ref{l3.2} hold true. Therefore, the following estimate holds obviously
   \bess
\sup_{(t,x,y)\in[0,T]\times\mathbb{R}^2}|v(t,x,y)|\leq T\sup_{(t,x,y)\in[0,T]\times\mathbb{R}^2}|f(t,x,y)|.
   \eess

Next we will use the rotation invariant of Laplacian operator and the estimates of
Lemma \ref{l3.2} to derive the desire results. In order to do that, we define
$S$ as an orthogonal transformation of $\mathbb{R}^2$: $Se_i=l_i$, $i=1,2$, where
$e_i$ is the standard basis in $\mathbb{R}^2$, $l_i$ is a unit vector in $\mathbb{R}^2$ and
$l_2$ is orthogonal to $l_1$. Set
   \bess
&&f(t,xe_1+ye_2)=f(t,x,y),\ v(t,xe_1+ye_2)=v(t,x,y),\ S(x,y)=xl_1+yl_2,\\
&& g(t,x,y)=f(t,S(x,y)),
w(t,x,y)=v(t,S(x,y)),
  \eess
then $w$ satisfies
   \bess
\partial_tw(t,x,y)=a(t)\Delta w(t,x,y)+g(t,x,y),
   \eess
where we used the rotation invariant of Laplacian operator.

It follows from Lemma \ref{l3.2} that
  \bess
&&\sup_{(t,y)\in[0,T]\times\mathbb{R}}\sup_{x_1\neq x_2}\frac{|D^2_{x}w(t,x_1,y)-D^2_{x}w(t,x_2,y)|}{|x_1-x_2|^\alpha}\\
&\leq&N(\alpha)\sup_{(t,y)\in(0,T)\times\mathbb{R}}\sup_{x_1\neq x_2}\frac{|g(t,x_1,y)-g(t,x_2,y)|}{|x_1-x_2|^\alpha}.
  \eess
Notice that
   \bess
D_x^2w(t,x,y)=D_{l_1}^2v(t,S(x,y))=D_{l_1}^2v(t,xl_1+yl_2),
   \eess
and using the fact that the solution $v$ of (\ref{2.6}) has continuous second-order derivatives w.r.t. $(x,y)$, we
have, for any unit vector $l\in\mathbb{R}^2$
    \bess
&&\sup_{(t,z)\in[0,T]\times\mathbb{R}^2}\sup_{\mu\neq \nu}\frac{|D^2_{l}v(t,\mu l+z)-D^2_{l}v(t,\nu l+z)|}{|\mu-\nu|^\alpha}\\
&\leq&N(\alpha)\sup_{(t,z)\in(0,T)\times\mathbb{R}^2}\sup_{\mu\neq \nu}\frac{|f(t,\mu l+z)-f(t,\nu l+z)|}{|\mu-\nu|^\alpha}.
  \eess
That is to say, we get
  \bess
\sup_{(t,z)\in[0,T]\times\mathbb{R}^2}[D^2_{l}v(t,z+l\cdot)]_{C^{\alpha}(\mathbb{R}^2)}\leq
N(\alpha)\sup_{(t,z)\in(0,T)\times\mathbb{R}^2}[D^2_{l}f(t,z+l\cdot)]_{C^{\alpha}(\mathbb{R}^2)}.
  \eess
In particular, if we choose $z=0$, we get the estimate
   \bess
\sup_{t\in[0,T]}[D_{ij}v(t,\cdot)]_{C^{\alpha}(\mathbb{R}^2)}\leq N_0(\alpha)\sup_{t\in[0,T]}[f(t,\cdot)]_{C^{\alpha}(\mathbb{R}^2)}.
  \eess
Since the Jacobian of $S(x,y)$ equals to $1$, then we have for any unit vector $l\in\mathbb{R}^2$
   \bess
\int_0^T\int_{\mathbb{R}^2}|D_l^2v(t,z)|^pdzdt\leq N_p
\int_0^T\int_{\mathbb{R}^2}|f(t,z)|^pdzdt.
   \eess
The proof is complete. $\Box$

\medskip

\noindent {\bf Acknowledgment} The first author was supported in part
by NSFC of China grants 11771123, 11531006, 11501577.

 \end{document}